\documentclass[aps,pre,reprint,groupedaddress]{revtex4-2}
\usepackage{amsmath,amssymb,bm}
\usepackage{graphicx}
\usepackage{hyperref}
\usepackage{color}
\usepackage{float}
\newcommand{\xq}[1]{\textcolor{red}{#1}}
\newcommand{\als}[1]{\textcolor{blue}{#1}}

\definecolor{ReviewBrown}{rgb}{0.55,0.27,0.07}
\definecolor{ReviewGreen}{rgb}{0.0,0.42,0.16}
\newcommand{\reviewplain}[1]{#1}
\DeclareRobustCommand{\reviewold}[1]{{\color{ReviewBrown}\let\xq\reviewplain\let\als\reviewplain #1}}
\DeclareRobustCommand{\reviewnew}[1]{{\color{ReviewGreen}\let\xq\reviewplain\let\als\reviewplain #1}}

\begin{document}
	
	\title{Delay-induced multirhythmicity in a 3-cell inhibitory network}
	
	
	\author{Xinxin Qie}
	\affiliation{School of Mathematics,
		South China University of Technology,
		Guangzhou 510640, China}
	\affiliation{Department of Information Engineering,
		University of Padova,
		35121 Padova, Italy}
	
	\author{Morten Gram Pedersen}
	\affiliation{Department of Information Engineering and Padova Neuroscience Center, University of Padova, 35121 Padova, Italy}
	
	\author{Jian Song}
	\affiliation{Department of Mathematics and Physics, School of Biomedical Engineering, Southern Medical University, Guangzhou 510515, China}
	
	\author{Shenquan Liu}
	\email{Contact author: mashqliu@scut.edu.cn}
	\affiliation{School of Mathematics, South China University of Technology, Guangzhou 510640, China}
	
	\author{Andrey L. Shilnikov}
	\email{Contact author: ashilnikov@gsu.edu}
	\affiliation{Neuroscience Institute and Department of Mathematics \& Statistics, Georgia State University, Atlanta, Georgia 30303, USA}
	\date{\today}
	
	\begin{abstract}
		We study the dynamics of a delayed three-cell central pattern generator (CPG) network using return-map analysis and a phase-reduced description derived from generalized FitzHugh–Nagumo oscillators. The phase model captures the emergence and stability of synchronous, pacemaker, and traveling-wave rhythms and provides a systematic, computationally efficient framework for characterizing the network dynamics beyond direct return-map calculations. We conduct extensive two-parameter analyses of both the full and phase-reduced models, revealing rich multi-rhythmic dynamics and transitions among coexisting rhythmic states. The close agreement between the two descriptions establishes phase reduction as an efficient and reliable framework for characterizing the dynamical organization of delayed neural oscillator networks.
	\end{abstract}
	
	\maketitle

	\section{Introduction}
	Central pattern generators (CPGs) are neural circuits that produce rhythmic motor patterns without rhythmic external input \cite{marder2001central,yuste2005cortex,kiehn2003physiological}. Interactions within biological CPGs are not instantaneous: finite synaptic transmission and signal-propagation times introduce delays into network coupling \cite{shafiei2020time,yu2018effect,venugopal2022connectionist}. Such delays can substantially reorganize the network dynamics, affecting synchronization, phase locking, and the coexistence of rhythmic states even in small networks of coupled relaxation oscillators \cite{qie2025synaptic,song2026coupling,campbell2012phase}.
	
	To analyze multistability in coupled CPG networks, Shilnikov and collaborators developed a return-map framework that provides a geometric description of phase-locked rhythms and transitions among them in small neuronal motifs \cite{wojcik2014key}. This approach has proved particularly effective for identifying coexisting rhythms, determining their basins of attraction, and characterizing the global organization of the network dynamics \cite{kelley20202,collens2020dynamics,wojcik2011order,schwabedal2016qualitative,pusuluri2020computational}. In our previous work \cite{qie2025synaptic}, we extended this framework to delayed inhibitory networks, enabling a systematic analysis of how synaptic delays affect rhythm generation, multistability, and transitions among rhythmic states.
	
	Phase reduction provides a complementary description of coupled-oscillator dynamics. In the weak-coupling regime, a high-dimensional network of oscillators can be reduced to a low-dimensional system governing their phase interactions \cite{ermentrout2010mathematical,kuramoto2005self,scully2024pairing,zhu2026dynamical}. Early studies of delay-coupled oscillators incorporated transmission delays explicitly into the phase arguments \cite{schuster1989mutual,yeung1999time,earl2003synchronization}. Subsequent work established that, under weak coupling, these delays can instead be represented as effective phase shifts in the interaction functions of the reduced equations \cite{crook1997role,bressloff1997phase,campbell2012phase,kurebayashi2022phase,novivcenko2012phase,kotani2012adjoint,ermentrout2010mathematical} as long as the delay is not too long \cite{Izhikevich1998}. 
	
	Previous studies have shown that, in CPG networks without delays, phase-difference models accurately reproduce the fixed-point organization and stability properties revealed by return maps \cite{zhu2026dynamical}. Whether this correspondence persists in delayed CPG networks, however, remains largely unexplored. Delays can generate more complex dynamical structures, including invariant curves, whose location and stability are not readily determined from return-map constructions alone. This limitation motivates a systematic comparison of return-map and phase-reduced descriptions of delayed CPG networks.
	
	In this work, we investigate a delayed three-cell inhibitory CPG network by combining return-map analysis with phase reduction of a network of generalized FitzHugh–Nagumo oscillators \cite{schwabedal2016qualitative,collens2020dynamics,pusuluri2020computational,qie2025synaptic}. Using phase-response curves (PRCs) and coupling functions \cite{ermentrout2010mathematical}, we derive a two-dimensional phase-difference system and systematically compare its predictions with the dynamics of the full delayed network.
	We show that the phase-reduced model not only reproduces the principal rhythmic structures revealed by the return maps but also provides a more detailed characterization of the phase space. In particular, it identifies stable synchronous states with extremely small basins of attraction that can be difficult to detect through direct return-map simulations. The reduced model also permits the location and stability of limit cycles in the phase-difference plane to be determined directly, thereby clarifying the dynamical origin of the invariant curves observed in the return maps.
	Using both the full and reduced models, we perform extensive two-parameter analyses of rhythm generation and transitions among coexisting states. The close agreement between the two descriptions establishes phase reduction as an accurate and computationally efficient framework for large-scale exploration of delayed CPG networks.

	\section{Methods}
	We consider three weakly coupled oscillators arranged in a symmetric networking ring, with a transmission delay in each coupling term. The state of the $j$-th oscillator, denoted by $X_j(t)\in\mathbb{R}^n$, evolves according to
	\begin{equation}
		\dot{X}_j
		=
		F(X_j)
		+
		g_{j-1,j}G_j^-
		+
		g_{j+1,j}G_j^+, j=1,2,3,
		\label{eq:network}
	\end{equation}
	where
	\begin{equation}
		\begin{split}
			G_j^-=G(X_j(t),X_{j-1}(t-\tau)),\\
			G_j^+=G(X_j(t),X_{j+1}(t-\tau)).
		\end{split}
	\end{equation}
	Here, $F:\mathbb{R}^n\to\mathbb{R}^n$ describes the intrinsic dynamics of each oscillator, $\tau\geq 0$ denotes the transmission delay, and $0<g_{j-1,j}\ll1$
	are small coupling parameters. The functions $G:\mathbb{R}^n\times\mathbb{R}^n\to\mathbb{R}^n$ specifies the interaction between a postsynaptic oscillator and a delayed presynaptic oscillator. 
	To identify activation events, we introduce a scalar observable
	\begin{equation}
		V_j(t)=\mathcal V(X_j(t)),
		\label{eq}
	\end{equation}
	where $\mathcal V:\mathbb{R}^n\to\mathbb{R}$ maps the state of each oscillator to the scalar component or output variable used to detect threshold crossings.
	
	Time delays are an intrinsic feature of neural signal transmission in biological systems, and their inclusion in CPG models is essential for accurately capturing signal propagation and enhancing biological plausibility. 
	Let $T_0$ denote the period of the uncoupled oscillator. Then, we can normalize the transmission delay $\tau$ by this intrinsic period as following: $\alpha=\tau/T_0$. 
	We will assume that $0 \le \alpha \le 1$.
	
	For a periodic phase-locked solution of the full network, let $T$ denote its common period. On a $T$-periodic network solution, $X_j(t-kT)=X_j(t)$ for every integer $k$.
	Thus, for $\tau=kT\geq0$, the delayed input coincides with the instantaneous input along that solution.
	This identity does not imply identical stability or transient dynamics in the delayed and undelayed networks.
	Noninteger-period delays generally modify the relative timing of the inputs and can reorganize the rhythmic states.
	\subsection{Return map}
	
	Because the full higher-dimensional network comprises three coupled oscillatory cells, its relative phase dynamics can be characterized by two independent phase-difference variables, $\Delta=\left(\Delta_{12}, \Delta_{13}\right)$. Here, $\Delta_{1 j}(j=2,3)$ denotes the phase lag of node $j$ relative to node $1$. Together, these variables constitute the maximal set of independent phase differences among the three oscillators.
	Following the procedure of Wojcik et al.~\cite{wojcik2011order}, we obtain the phase lags from threshold-crossing times of the membrane potentials. Specifically, we choose node 1 as the reference oscillator and record the sequence of times $t_1^{(k)}$ at which $V_1(t)$ crosses the threshold $V_1=0$ from below, with each crossing corresponding to the $k_{th}$ activation of node 1. For each reference event, we then identify the subsequent activation times $t_j^{(k)}$, $j=2,3$, of nodes 2 and 3. The time lags of these events relative to the reference event are normalized by the corresponding reference period $T_1^{(k)}=t_1^{(k+1)}-t_1^{(k)}$. This procedure yields a sequence of phase-lag pairs along the forward trajectory, denoted by $\Delta^{(k)}=\left(\Delta_{12}^{(k)},\Delta_{13}^{(k)}\right)$, where $\Delta_{12}^{(k)}=(t_2^{(k)}-t_1^{(k)})/T_1^{(k)},\Delta_{13}^{(k)}=(t_3^{(k)}-t_1^{(k)})/T_1^{(k)}.$
	The cycle-to-cycle evolution of the phase lags is governed by the Poincaré map $\Pi: P_{k}\mapsto P_{k+1}$, where $P_k=\Delta^{(k)}$. This map advances the phase-lag vector from one activation cycle of the reference node to the next. Its dynamics determine whether trajectories converge to stable fixed points (FPs), corresponding to phase-locked rhythms with constant phase lags, or exhibit more complex behavior, such as motion on an invariant closed curve in the phase-lag plane. Stable FPs represent asymptotically stable phase-locked rhythms, whereas the coexistence of multiple attractors indicates network multistability.
	
	When represented on the unit square with both coordinates defined modulo 1, the two-dimensional phase-lag map preserves the topology of a torus, with opposite edges identified. At the pairwise level, the basic phase relationships include $\Delta=0$ for cells in-phase and $\Delta=0.5$ for cells in anti-phase. To identify the attractors and their basins of attraction, we initialize the network on a grid of phase-lag values (e.g., $20\times20$
	) and compute sufficiently long sequences of network activity for the transient dynamics to decay. Fixed points with constant $(\Delta_{12}^*, \Delta_{13}^*)$ lags correspond to stable phase-locked rhythms. Initial conditions assigned the same color converge to the same attractor and therefore belong to the same basin of attraction. Distinct rhythmic states are represented by fixed points at different locations in the phase-difference plane.

	We classify the observed rhythms according to the locations and stability of the fixed points (FPs) and the geometry of the asymptotic trajectories in the phase plane. Traveling-wave (TW) rhythms are associated with stable FPs near $(\Delta_{12},\Delta_{13})\approx(1/3,2/3)$ and its symmetry-related counterpart. Pacemaker (PM) rhythms correspond to stable FPs near $(0,1/2)$, $(1/2,1/2)$, and $(1/2,0)$, whereas synchronous (SY) rhythms are identified by convergence to $(\Delta_{12},\Delta_{13})\approx(0,0)$. Trajectories that converge to a closed invariant curve in the phase plane are associated with the phenomenon of phase jitter -- slow oscillations centered around a phase-locked state, whereas trajectories exhibiting persistent phase drift are classified as phase-slipping (PS) dynamics.
	\begin{figure}[t!]
		\includegraphics[width=0.5\textwidth]{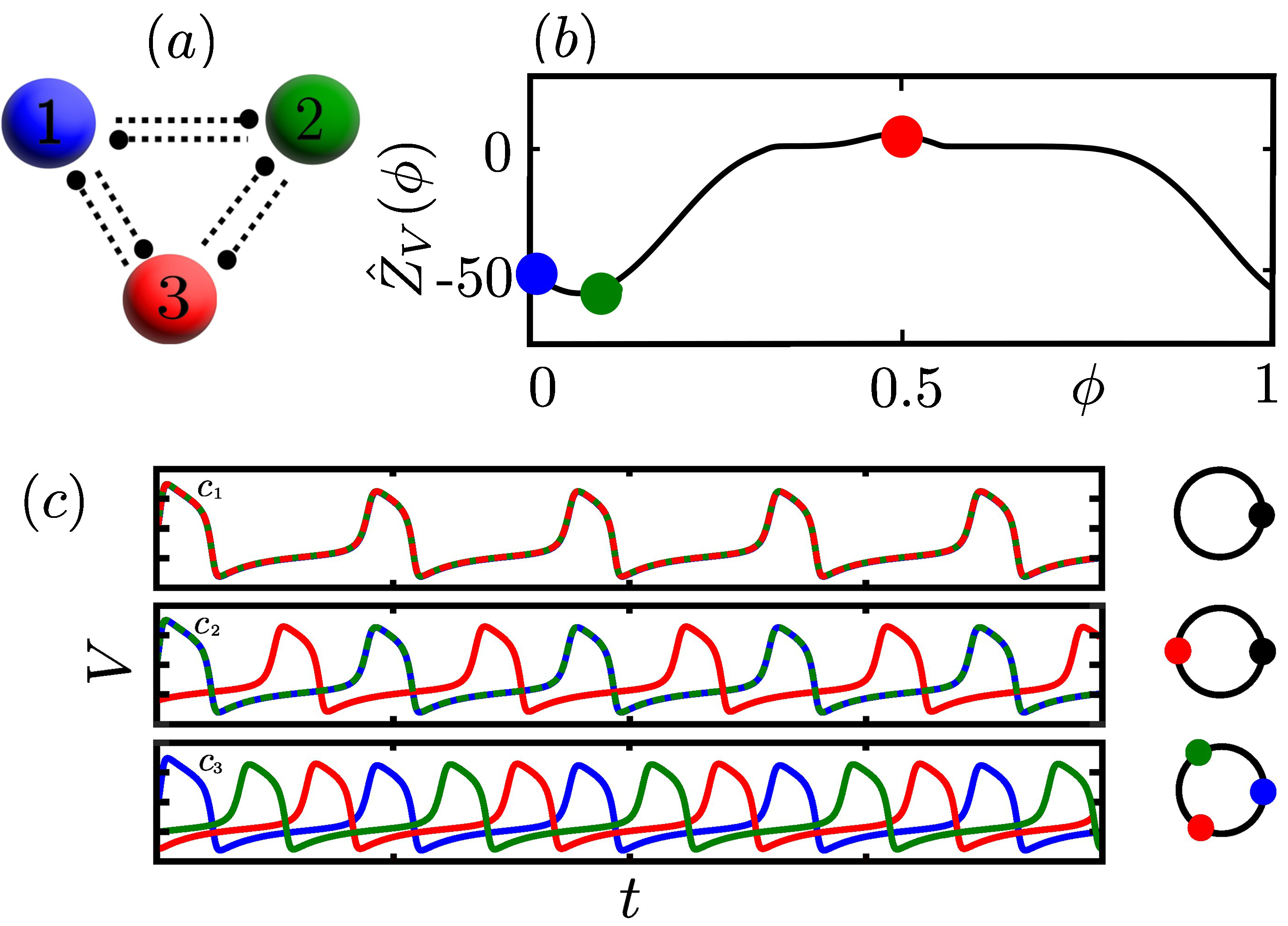}
		\caption{\label{PRC}(a) A neural motif of three oscillatory cells coupled by delayed inhibitory synapses. (b) Voltage component of the phase-response curve (PRC), $\hat{Z}_V(\phi)$, for $I=0.6$ and $\varepsilon=0.3$. Cell 1 (blue) is taken as the reference, with $\phi_1=0$, while the colored circles illustrate the phase positions of the three cells. Repeated delayed inhibitory interactions modify
			the relative phases according to the PRC.
			(c) Three samples of typical ``voltage'' traces generated by the 3-cell network: (c1)  a synchronous rhythm with $\Delta_{12}=\Delta_{13}=0$; (c2) a pacemaker rhythm (PM) where one red cell oscillates in anti-phase ($\Delta_{13}=1/2$); with two synchronous ($\Delta_{12}=0$) blue and green cells; (c3)  a clockwise traveling wave (TW) with the following phase-locked states $\Delta_{12}\simeq 1/3$ and $ \Delta_{13}\simeq 2/3$.} 
\end{figure}
\begin{figure}[t!]
	\centering 
	\includegraphics[width=0.45\textwidth]{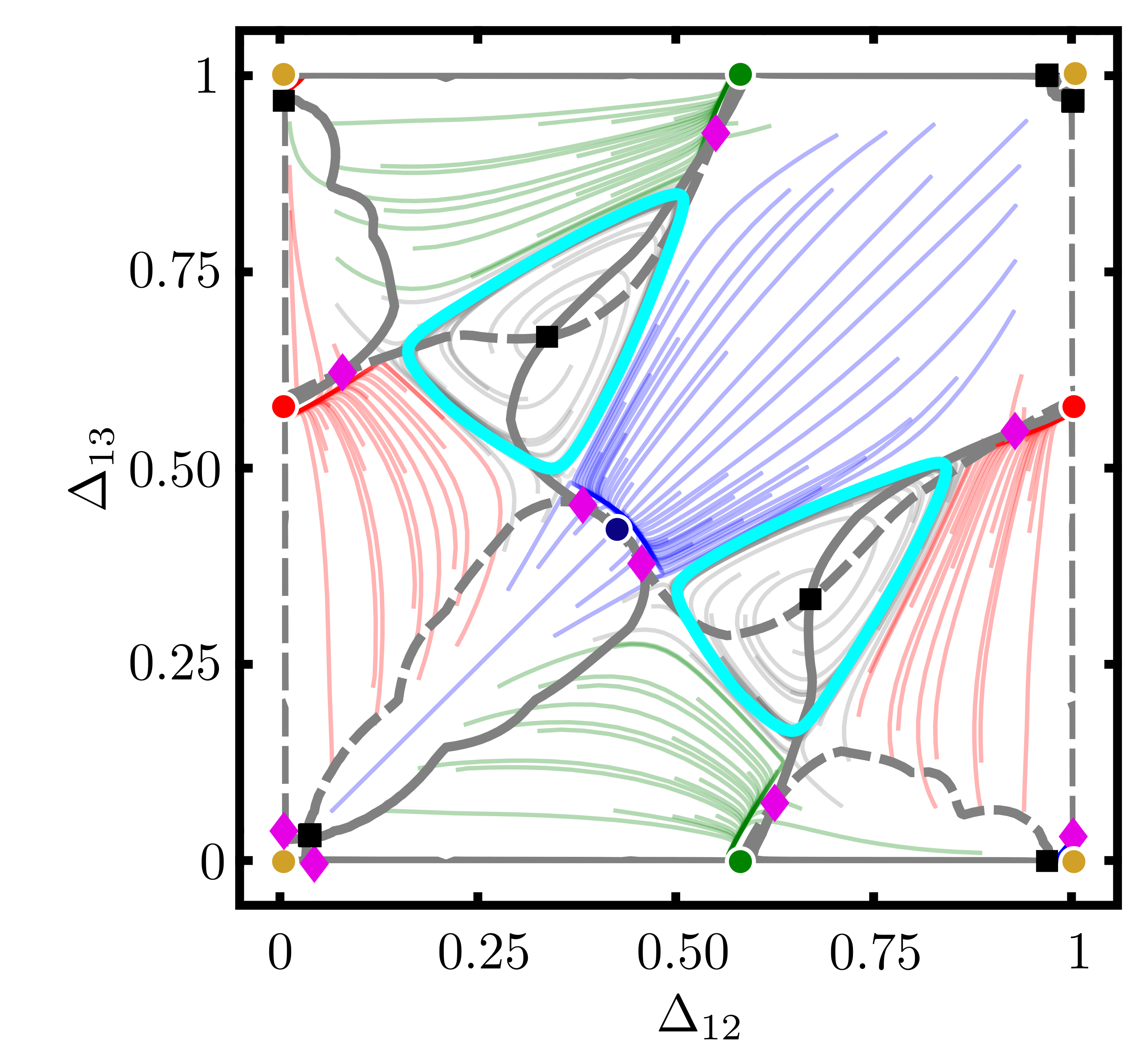}
	\caption{\label{compare} Comparison between the phase-difference dynamics obtained from the full model and the phase
		reduced model. The thick solid and dashed curves represent the nullclines obtained from the phase--
		difference model (\ref{eq:phase_difference_system}). 			The low-opacity trajectories are computed from the return map and illustrate the asymptotic evolution toward the attractors. The stable fixed points (FPs) located at crossings of both nullclines almost perfectly coincide with the attractors obtained from the return map, demonstrating that the phase model can accurately determine the locations of the stable states in the phase space. The red, green and blue solid circles represent the PMs. The yellow dots stand for the stable SY at the origin, while the black square $\blacksquare$ and purple diamond $\blacklozenge$ symbols represent unstable FPs and saddles, and the cyan curves represent the stable invariant curves (ICs).  
		Parameters: $I=0.6$, $\varepsilon=0.3$, $\alpha=0.04$ and $g=0.001$.  
	}
\end{figure}

\subsection{Standard phase reduction}
We assume that each uncoupled oscillator possesses an exponentially asymptotically stable limit cycle $\widehat X(t)$ with period $T_0$. Let $\widehat Z(t)$ denote the corresponding $T_0$-periodic infinitesimal phase-response curve (PRC), normalized according to
\begin{equation}
	\widehat Z(t)\cdot F(\widehat X(t))=\frac{1}{T_0}.
	\label{eq:prc_normalization}
\end{equation}
For sufficiently weak coupling, the delayed three-cell network introduced above can be reduced to a phase model~\cite{ermentrout2010mathematical}. 

The interaction function $H$ is defined by
\begin{equation}
	\label{eq:H}
	H(\phi)=\frac{1}{T_0}\int_0^{T_0}
	\widehat Z(s)\cdot
	G\bigl(\widehat X(s),\widehat X(s+\phi T_0)\bigr)\,ds.
\end{equation}

Because the phase variables are defined modulo $1$ and $\widehat X$ is $T_0$-periodic, the interaction function $H$ is $1$-periodic.
The PRC determines the phase sensitivity of each oscillator to weak coupling perturbations. 
For an oscillator at a given phase, as illustrated by the colored circles in Fig.~\ref{PRC}(b), each delayed coupling input produces a small phase correction whose magnitude and direction depend on the phase at which the perturbation is received.
Through repeated interactions, these phase corrections continuously modify the relative phases between the oscillators.
Depending on the initial phase and coupling function, the network can exhibit different stable phase-locked rhythms, with representative examples shown in Fig.~\ref{PRC}(c).


Denoting the normalized phase of oscillator $i$ by $\varphi_i\in[0,1)$, to leading order in the weak-coupling limit,
we obtain \cite{ermentrout2010mathematical} 
\begin{equation}
	\frac{d\varphi_i}{dt}
	=\Omega+\sum_{\substack{j=1\\j\neq i}}^{3}g_{ji}
	H(\varphi_j-\varphi_i-\alpha),
	\qquad i=1,2,3,
	\label{eq:phase_model}
\end{equation}
where $\Omega=1/T_0$ and $g_{ji}$ 
is the 
coupling strength from oscillator $j$ to oscillator $i$.
Since one uncoupled oscillation corresponds to one unit of phase,
the transmission delay $\tau$ 
corresponds to the phase delay
\[
\alpha=\Omega\tau=\frac{\tau}{T_0}\pmod{1},
\]
where $T_0$ is the period of the uncoupled oscillator and
$\Omega=1/T_0$. 
When $\alpha$ is small compared to $1/(\max g_{ij})$ it enters the reduced equations as a
fixed phase shift in the interaction function \cite{Izhikevich1998}.
The resulting system~(6) consists of three phase equations; eliminating the common phase yields a two-dimensional system governing the two independent phase differences.

Choose phase zero at the reference threshold crossing, so that $\varphi_1(t_1)=0$ modulo one. For uncoupled oscillators the time-to-phase conversion is exact with period $T_0$. Under weak coupling it holds to leading order for transverse crossings and bounded inter-event intervals:
\begin{equation}
	\varphi_j(t_1)=1-\frac{t_j-t_1}{T_0}+O(g)\pmod 1.
\end{equation}
Hence, defining the reduced phase differences by
\begin{equation}
	\label{eq:phase_differences}
	\begin{aligned}
		\Delta_{1j}&=\varphi_1-\varphi_j\\
		&\equiv\frac{t_j-t_1}{T_0}+O(g)\pmod 1,
		\qquad j=2,3,
	\end{aligned}
\end{equation}
gives the same sign convention as the event-based lags. When $T_1^{(k)}=T_0+O(g)$, normalization by the measured reference interval instead of $T_0$ also changes these lags only by $O(g)$. The reduced phases and the measured event lags are thus leading-order counterparts, rather than exactly identical coordinates.
These variables evolve at rates of order $g$, and hence on a slow time scale of order $g^{-1}$. To leading order in $g$, they satisfy
\begin{equation}
	\label{eq:phase_difference_system}
	\begin{aligned}
		\dot{\Delta}_{12}
		&=
		g_{21}H(-\Delta_{12}-\alpha)
		+
		g_{31}H(-\Delta_{13}-\alpha)
		\\
		&\qquad
		-
		g_{32}H(\Delta_{12}-\Delta_{13}-\alpha)
		-
		g_{12}H(\Delta_{12}-\alpha)
		,
		\\[1ex]
		\dot{\Delta}_{13}
		&=
		g_{21}H(-\Delta_{12}-\alpha)
		+
		g_{31}H(-\Delta_{13}-\alpha)
		\\
		&\qquad
		-
		g_{23}H(\Delta_{13}-\Delta_{12}-\alpha)
		-
		g_{13}H(\Delta_{13}-\alpha)
		.
	\end{aligned}
\end{equation}

The fixed points (FPs), given by the intersections of the nullclines of Eq.~\eqref{eq:phase_difference_system}, correspond to phase-locked solutions of the original oscillator network. They are located at the intersections of the two nullclines,  which are defined by the conditions: $\dot{\Delta}_{12}=0$ and $\dot{\Delta}_{13}=0$, and are shown as the red and blue curves, respectively, in Fig.~\ref{compare}. Linearization of Eq.~\eqref{eq:phase_difference_system} about a FP $\Delta^*$ yields the Jacobian matrix $J$. The eigenvalues of $J$ determine the local stability of the corresponding phase-locked rhythm.

\subsection{Numerics}

Simulations of the full network equations were performed with a fixed time step $dt=0.01$. The instantaneous state variables were advanced using a classical fourth-order Runge--Kutta scheme, while the delayed presynaptic terms were evaluated from a discrete history buffer. For $t\in[-\tau,0)$, the delayed presynaptic variables were initialized to zero, whereas the oscillator states at $t=0$ were specified by the chosen initial conditions. Within each integration step, the delayed input obtained from the history buffer was held fixed across the four Runge--Kutta stages, without additional interpolation. Step-size convergence was checked by repeating the simulations with a smaller time step, with no appreciable changes in the reported dynamical behavior.
The resulting trajectories are processed using Python scripts for visualization and return-map construction. The nullclines of the phase-difference system obtained using MATLAB were double-checked with calculations using XPPAUT~\cite{ermentrout02xppaut}. Through systematic parameter variation, we track changes in the geometry of the return map and identify bifurcations and transitions among the network rhythms. This analysis reveals dynamical structures that are not readily inferred from voltage traces alone and characterizes the parameter dependence, stability, and multistability of the network-generated rhythms.

\section{Results}
\subsection{Generalized FitzHugh--Nagumo network}
To compare and validate the return-map and phase-reduction approaches, we use a generalized FitzHugh--Nagumo (FHN) model as a concrete realization of each neuronal oscillator. The FHN model is a canonical low-dimensional reduction of excitable membrane dynamics and can be viewed as a biologically motivated extension of the van der Pol relaxation oscillator. It retains the essential interaction between fast voltage activation and slow recovery while omitting the detailed ionic conductances of Hodgkin--Huxley-type models. The generalized FHN model used here incorporates an explicit slow--fast time-scale separation and a sigmoidal recovery nullcline, allowing it to reproduce relaxation oscillations and the envelope-level alternation between active and quiescent phases observed in neuronal bursters, both in isolation and under synaptic perturbations~\cite{collens2020dynamics}:
\begin{equation}
	\begin{aligned}
		\dot V_i
		&=
		V_i-V_i^3-h_i+I
		+
		\sum_{j\ne i}g_{ji}G_{ji}(V_j,V_i),
		\\
		\dot h_i
		&=
		\varepsilon
		\left[
		\frac{1}{1+e^{-k(V_i-V_0)}}-h_i
		\right],
		\qquad i=1,2,3.
	\end{aligned}
	\label{eq:fhn_network}
\end{equation}
Here, the state of the \(i\)th neuron is described by the fast membrane-voltage variable \(V_i\) and some  slower recovery or voltage-dependent gating  \(h_i\)-variable. The small parameter \(0<\varepsilon <1 \) controls the separation between the fast voltage and slow recovery time scales, whereas \(I\) is an applied current. The parameter \(V_0\) determines the midpoint of the sigmoidal \(h\)-nullcline, and \(k\) controls its steepness. The cubic \(V\)-nullcline is determined by the voltage equation and depends on \(I\) and the synaptic input. Unless stated otherwise, we use \(k=10\) and \(V_0=0\). We also set the synaptic activation threshold to \(V_{\mathrm{th}}=0\). The time-scale parameter \(\varepsilon\) should not be confused with the weak-coupling parameter $g$
introduced in the general network formulation.

A presynaptic neuron activates its outgoing synapse when its voltage exceeds \(V_{\mathrm{th}}\). Synaptic interactions are modeled using the fast threshold modulation (FTM) formalism,
\begin{equation}
	\begin{aligned}
		G_{ji}(V_j,V_i)
		&=
		(V_{\mathrm{rev}}-V_i)\Gamma(V_j),
		\\
		\Gamma(V_j)
		&=
		\frac{1}
		{1+e^{-100\left(V_j(t-\tau)-V_{\mathrm{th}}\right)}}.
	\end{aligned}
	\label{eq:ftm_synapse}
\end{equation}
Here, \(\tau\) is the transmission delay and is treated as a bifurcation parameter. The steep sigmoidal function \(\Gamma\) approximates a binary synaptic switch: it is close to \(1\) when the delayed presynaptic voltage exceeds \(V_{\mathrm{th}}\) and close to \(0\) when it lies below the threshold. Thus, the FTM formalism distinguishes between the active (``on'') and inactive (``off'') states of the synaptic gate.

The nonnegative weight \(g_{ji}\) specifies the strength of the directed synapse from neuron \(j\) to neuron \(i\). The sign of the synaptic current is determined by the driving force \(V_{\mathrm{rev}}-V_i\). A reversal potential below the relevant voltage range, such as \(V_{\mathrm{rev}}=-1.5\), produces inhibitory coupling, whereas a reversal potential above that range, such as \(V_{\mathrm{rev}}=1.5\), produces excitatory coupling. In the present study, we focus on inhibitory interactions.

It is well established~\cite{qie2025synaptic} that delayed three-node CPG networks can generate three canonical classes of rhythms: a synchronous state (SY; yellow point), three pacemaker states (PM; green, red, and blue points), and two traveling-wave states (TW; purple and black points), as shown in Figs.~\ref{compare} and~\ref{tau0.08}. Cyan lines and cyan circles represent phase-slipping (PS) and limit-cycle (LC) dynamics, respectively, with the corresponding phase trajectories shown in gray. Purple diamonds denote saddle FPs, whereas black squares denote unstable FPs.

We compare the phase-space structures obtained from the return map and the phase model in the $(\Delta_{12},\Delta_{13})$--phase-lag map. The two approaches show good agreement across several parameter sets, and one representative example is shown in Fig.~\ref{compare}.
The locations and stability of the FPs inferred from the return map agree closely with those determined from the nullcline intersections for locating FPs, while their stability is determined by local linearization of the phase model. Direct integration of the phase-difference system also identifies stable LCs corresponding to invariant closed curves. In some parameter regimes, the basin of attraction of the synchronous state at the origin is extremely small, making the stable yellow SY FP difficult to detect from return-map trajectories alone. The phase model resolves this difficulty by determining the stability of the corresponding FP directly from the eigenvalues of the linearized phase-difference system, independently of its basin size.

Using both return-map analysis and phase reduction, we systematically investigate the dynamics of the delayed three-neuron network across different parameter regimes. The return-map framework is used to analyze the effects of coupling strength and transmission delay directly in the full network. The corresponding phase model provides the precise locations and stability of saddle and unstable FPs and identifies the LCs associated with invariant curves in the return map. We also use the reduced model to explore the $(\varepsilon,I)$ parameter plane, for which its lower dimensionality permits denser parameter sweeps and a more detailed classification of the resulting rhythms at substantially lower computational cost. Together, the two approaches reveal a rich organization of coexisting rhythmic attractors and the bifurcations that mediate transitions among them.

\begin{figure}[t!]
	\centering
	\includegraphics[width=0.5\textwidth]{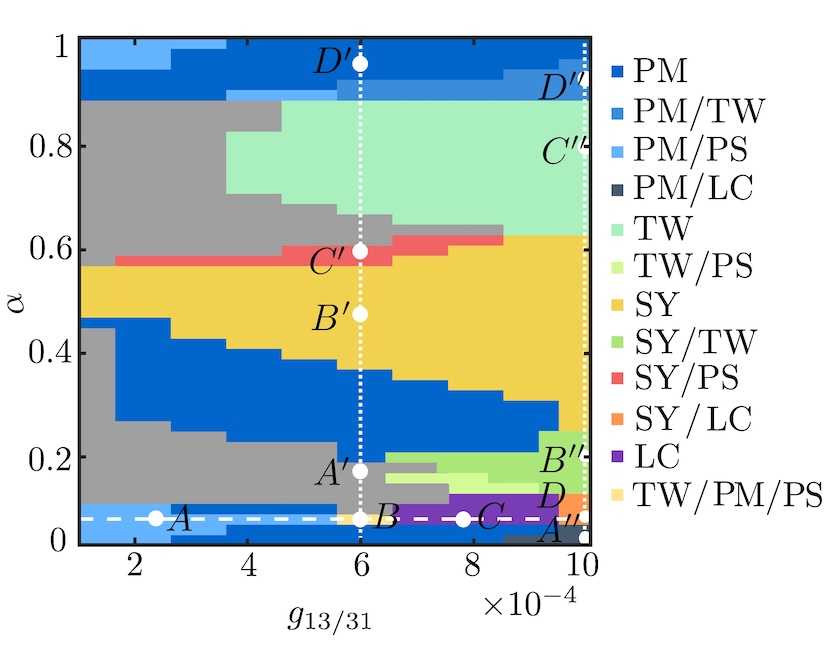}
	\caption{\label{g13_31}
		Color-coded dynamical regimes in the (\(g_{13/31},\alpha)\) parameter plane. The colors denote distinct rhythmic states and their coexistence, including PM, TW, SY, PS, and LC dynamics. Points \(A\)--\(D\), \(A'\)--\(D'\), and \(A''\)--\(D''\) indicate the parameter values used for the return maps in Figs.~\ref{tau0.08}, \ref{g0.0006}, and \ref{g0.001}, respectively.}
\end{figure}

\subsection{Return-map analysis of asymmetric inhibitory networks}
\begin{figure}[t!]
	\centering
	\includegraphics[width=0.5\textwidth]{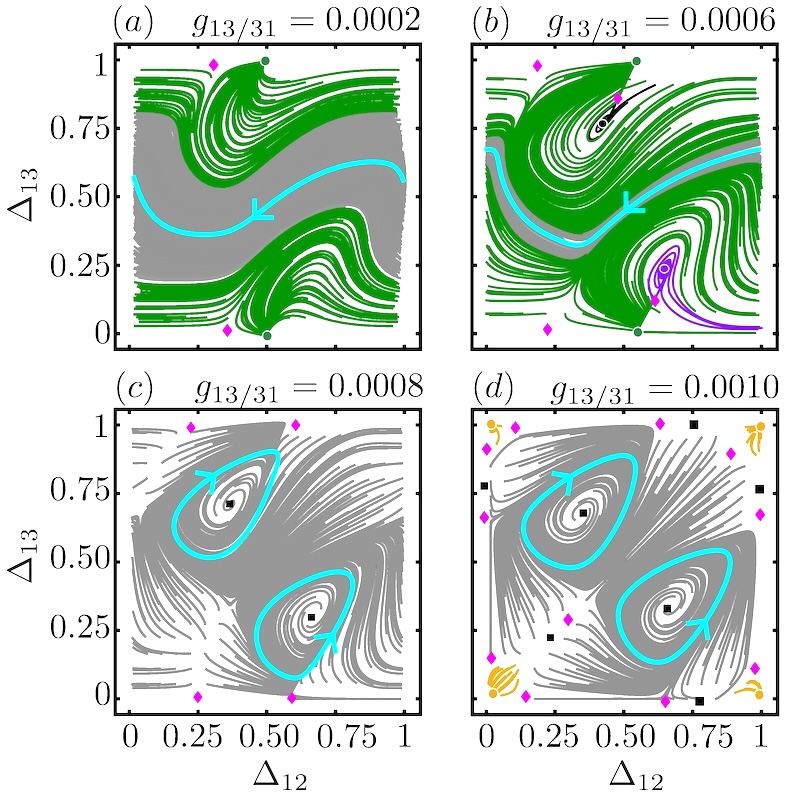}
	\caption{\label{tau0.08}
		Poincar\'e return maps for the parameter values \(A\)--\(D\) in Fig.~\ref{g13_31}. At weak reciprocal coupling, \(g_{13/31}=0.0002\), the green PM coexists with PS dynamics, shown in gray, in panel~(a). A representative PS trajectory is highlighted in cyan, with arrows indicating its direction. As \(g_{13/31}\) increases, the clockwise black and counter clockwise purple stable TWs emerge with initially small basins of attraction in panel~(b). With a further increase in coupling, the stable PM and TW FPs are no longer observed, while two stable LCs remain in panel~(c). When \(g_{13/31}=g_{ij}=0.001\), restoring coupling symmetry, a stable SY state with a small basin is observed in panel~(d). The remaining notation and color coding follow Fig.~\ref{compare}. The yellow dots stand for the stable SY at the origin, while the black square $\blacksquare$ and purple diamond $\blacklozenge$ symbols represent unstable FPs and saddles, and the cyan curves represent the stable invariant curves (ICs). 	 Parameters: \(\alpha=0.08\); all synaptic weights are fixed at \(g_{ij}=0.001\)
		, except for \(g_{13}=g_{31}=g_{13/31}=0.0002,0.0006,0.0008,\) and \(0.001\) in panels~(a)--(d), respectively.}
\end{figure}

We investigate how asymmetry in the inhibitory coupling affects the rhythmic dynamics of the delayed three-neuron network. Specifically, the reciprocal coupling pair \(g_{13}=g_{31}\equiv g_{13/31}\) is varied relative to the remaining synaptic weights. Using Poincar\'e return maps together with a systematic parameter sweep, we classify the multirhythmic dynamics in the two-parameter, \((g_{13/31},\alpha)\), plane.

Figure~\ref{g13_31} shows the global dynamical organization of the delayed three-neuron motif in the \((g_{13/31},\,\alpha)\)-parameter plane. The distribution of dynamical regimes is approximately reflected about \(\alpha=0.5\), although this symmetry is not exact. Because \(\alpha\) and \(1-\alpha\) correspond to phase shifts of opposite orientation modulo one period, the approximate reflection is consistent with the exchange symmetries of the ring architecture. Its imperfect character reflects the detailed shape of the interaction function and the nonlinear dynamics of the full network.

Pacemaker rhythms dominate for both short and nearly full-period delays, particularly for \(\alpha\in(0,0.1)\) and \(\alpha\in(0.9,1)\). In contrast, PS dynamics are distributed approximately symmetrically about \(\alpha=0.5\). Representative parameter regions and their characteristic rhythms are examined below. Unless stated otherwise, the intrinsic parameters are fixed at \(I=0.6\) and \(\varepsilon=0.3\) for all three neurons.

\begin{figure}[t!]
	\centering
	\includegraphics[width=0.5\textwidth]{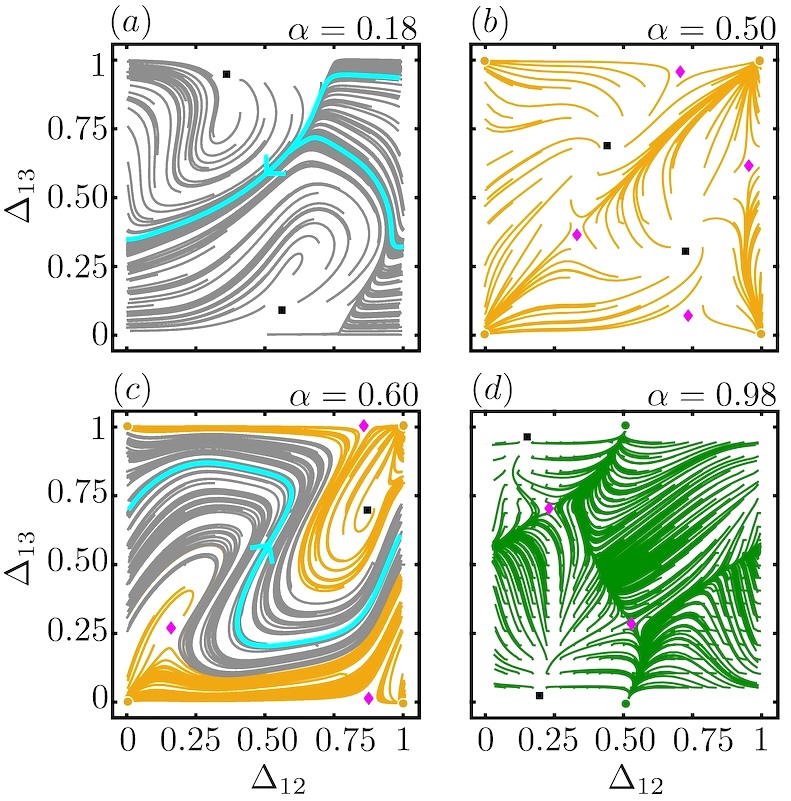}
	\caption{\label{g0.0006}
		Poincar\'e return maps for the parameter values \(A'\)--\(D'\) in Fig.~\ref{g13_31}. At \(\alpha=0.18\), PS dynamics dominate in panel~(a). A stable SY state -- the yellow dot near the origin appears as the delay increases in panel~(b) and coexists with PS dynamics in panel~(c). For a delay close to one oscillation period, a green PM appears and becomes the dominant rhythm in panel~(d). The notation and color coding follow Fig.~\ref{tau0.08}. Parameters: \(g_{13/31}=0.0006\) and \(\alpha=0.18,0.50,0.60,\) and \(0.98\) in panels~(a)--(d), respectively.}
\end{figure}

Figure~\ref{tau0.08} illustrates the Poincar\'e return maps for a short synaptic delay, \(\alpha=0.08\), as the reciprocal coupling strength \(g_{13/31}\) is varied. For weak coupling, see Fig.~\ref{tau0.08}(a), a stable green PM coexists with PS dynamics. Unlike a phase-locked attractor, the PS state does not converge to a single FP or invariant closed curve. Instead, trajectories initialized from different points follow distinct, non-convergent drifting paths on the phase torus, indicating persistent phase drift without synchronization or stable phase locking.

As \(g_{13/31}\) increases, two stable traveling wave FPs or TWs located symmetrically near \((1/3,2/3)\) and \((2/3,1/3)\) in the map and  shown in black and purple are observed in Fig.~\ref{tau0.08}(b). Although stable, these TWs initially possess small basins of attraction and are therefore strongly sensitive to the initial phase configuration. With a further increase in coupling, the PM and TW FPs disappear, and the dynamics become dominated by two stable LCs surrounding the former TW regions [Fig.~\ref{tau0.08}(c)]. When \(g_{13/31}=0.001\), the coupling symmetry is restored. A stable SY state is then observed with a very small basin of attraction, leaving the phase space predominantly occupied by the two LCs [Fig.~\ref{tau0.08}(d)].

Figure~\ref{g0.0006} shows the evolution of the Poincar\'e return maps for the asymmetric motif with \(g_{13/31}=0.0006\) as the synaptic delay increases. For relatively short delays, the network dynamics are dominated by PS [Fig.~\ref{g0.0006}(a)]. As the delay increases, a stable SY state is observed in Fig.~\ref{g0.0006}(b) and subsequently coexists with PS dynamics in Fig.~\ref{g0.0006}(c). Finally, when the delay approaches one full oscillation period, a green PM emerges and becomes the dominant rhythm [Fig.~\ref{g0.0006}(d)].

\begin{figure}[t!]
	\centering
	\includegraphics[width=0.5\textwidth]{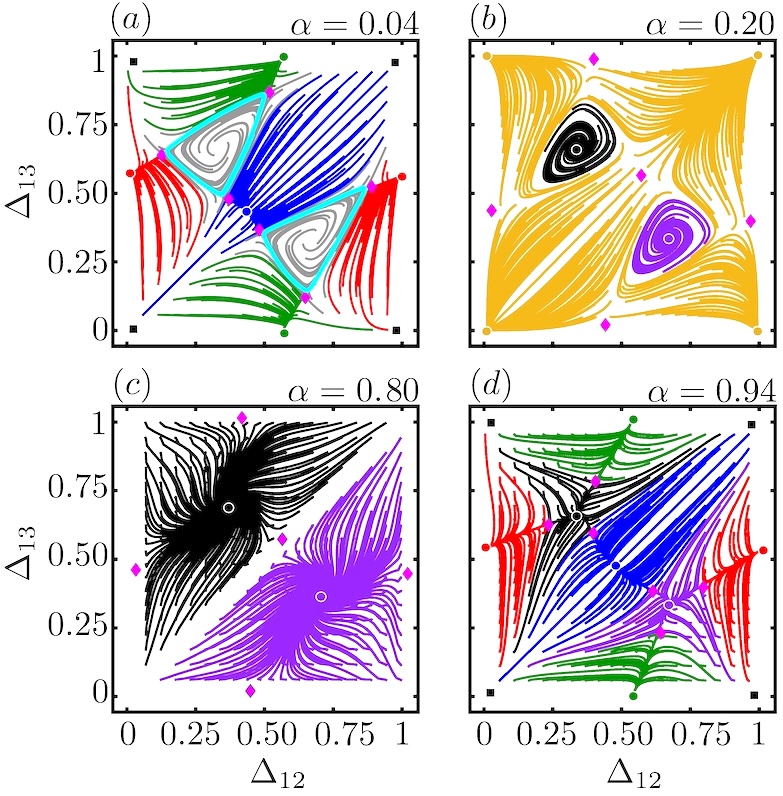}
	\caption{\label{g0.001}
		Poincar\'e return maps for the parameter values \(A''\)--\(D''\) in Fig.~\ref{g13_31}. For short delays, two stable LCs coexist with three stable PMs in panel~(a). As the delay increases, the yellow SY FP and the black and purple TWs emerge in panel~(b). The TWs become the dominant rhythms for longer delays in panel~(c), whereas the PMs reemerge as the delay approaches one full oscillation period in panel~(d). The notation and color coding follow Figs.~\ref{compare} and~\ref{tau0.08}. The yellow dots stand for the stable SY at the origin, while the black square $\blacksquare$ and purple diamond $\blacklozenge$ symbols represent unstable FPs and saddles, and the cyan curves represent the stable invariant curves (ICs).  Parameters: all synaptic weights satisfy \(g_{ji}=0.001\) and \(\alpha=0.04,0.20,0.80,\) and \(0.94\) in panels~(a)--(d), respectively.}
\end{figure}

\begin{figure}[t!]
	\includegraphics[width=0.5\textwidth]{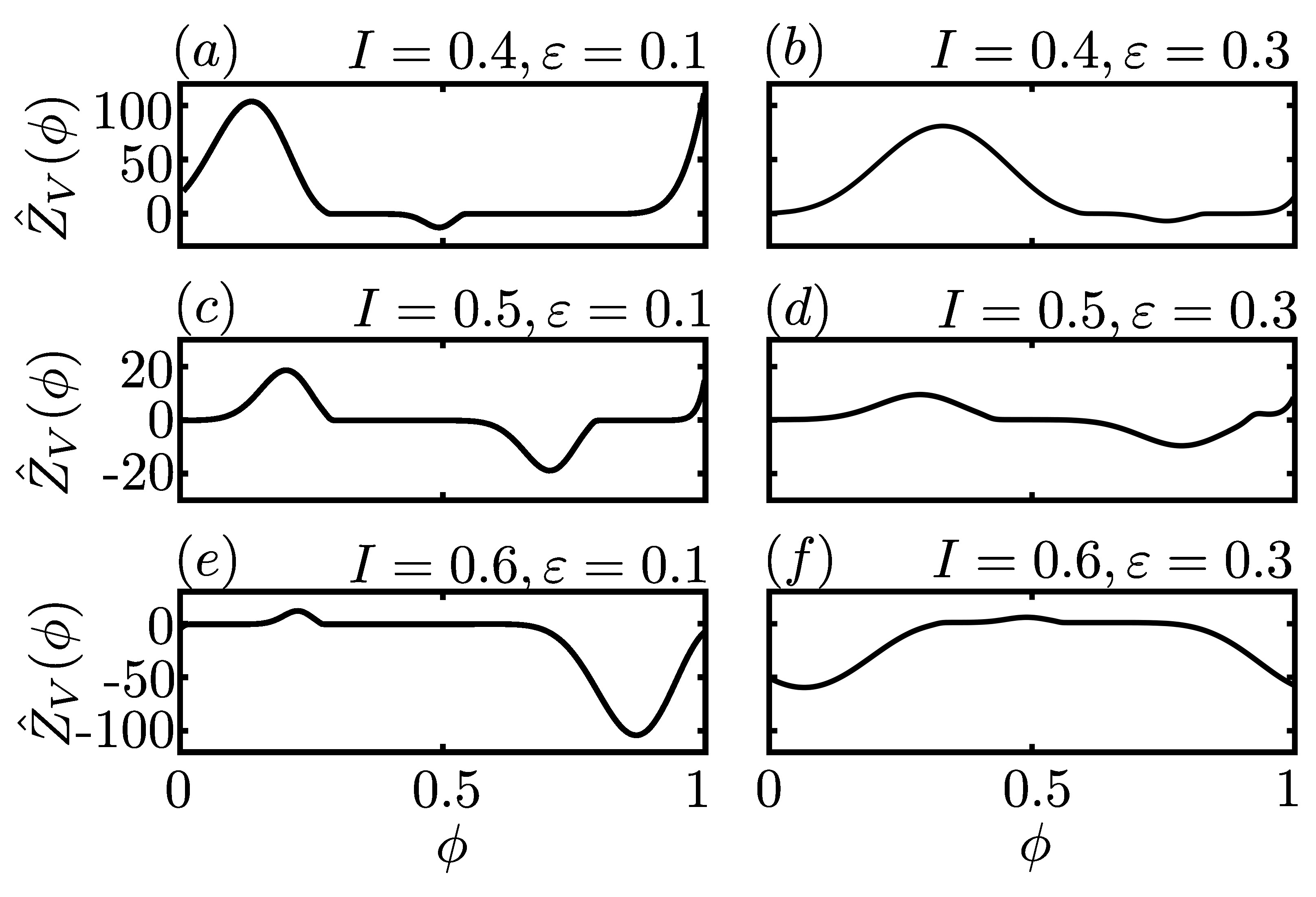}
	\caption{\label{PRCs}Phase-response curves (PRCs) $\hat{Z}_V (\phi)$ of the uncoupled FHN oscillator for different values of the intrinsic parameters $I$ and $\varepsilon$. For all parameter combinations, the limit-cycle orbit is computed starting from the same initial condition. The PRCs were computed from the adjoint equation along the uncoupled limit cycle and normalized according to Eq.~\eqref{eq:prc_normalization}. The rows correspond to $I=0.4$, $0.5$, and $0.6$, while the columns correspond to $\varepsilon=0.1$ and $0.3$.}
\end{figure}

\begin{figure*}[tb!]
	\centering
	\includegraphics[width=0.9\textwidth]{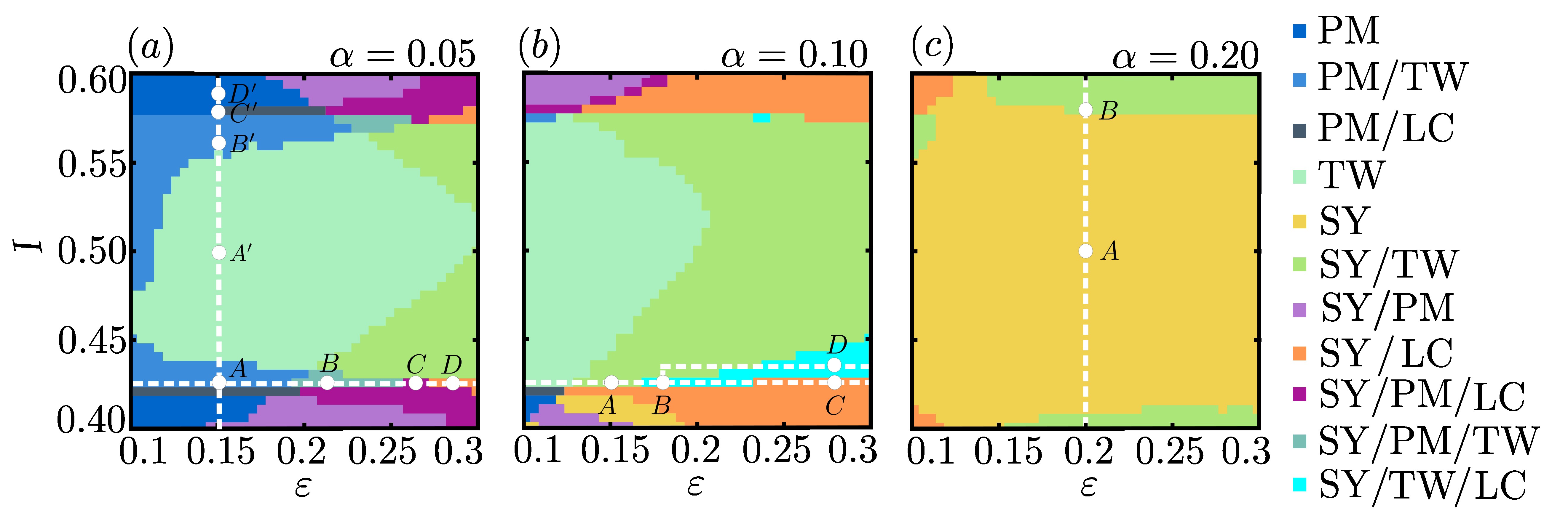}
	\caption{\label{3d}
		Two-parameter \((\varepsilon,I)\) sweeps at \(\alpha=0.05\), \(0.10\), and \(0.20\)
		in panels~(a)--(c), resp. The diagrams show the resulting rhythmic regimes in the phase difference model, including TW, PM, and SY states and their coexistence.}
\end{figure*}

Figure~\ref{g0.001} illustrates multistability in the symmetric motif with \(g_{ij}=0.001\).
For short delays, three PM states coexist with two stable LCs [Fig.~\ref{g0.001}(a)]. In the symmetric phase model, the SY equilibrium exists for every $\alpha$-value because all three phase velocities coincide when the phases coincide. Thus an observed onset of stable synchrony must be distinguished from creation of the equilibrium. As the delay increases, stable SY and TW states are observed [Fig.~\ref{g0.001}(b)]. For longer delays, the TWs dominate the phase space [Fig.~\ref{g0.001}(c)]. As the delay approaches one oscillation period, the PMs reappear and again become dominant [Fig.~\ref{g0.001}(d)]. This sequence is qualitatively consistent with the organization reported for \(g_{\mathrm{all}}=0.0005\) in Ref.~\cite{qie2025synaptic}, indicating that the principal polyrhythmic structure persists as the coupling strength is increased.

\subsection{Weak-coupling analysis via phase reduction}
Although the return-map approach provides a detailed description of the full network dynamics, its computational cost increases substantially in the weak-coupling regime. As the coupling strength decreases, trajectories approach their asymptotic states on a slow time scale of order \(g^{-1}\). Accurate classification of the attractors and their basins therefore requires long integration times at every point in a parameter sweep, making direct return-map construction inefficient for large-scale analyses.

\begin{figure}[b!]
	\centering
	\includegraphics[width=0.5\textwidth]{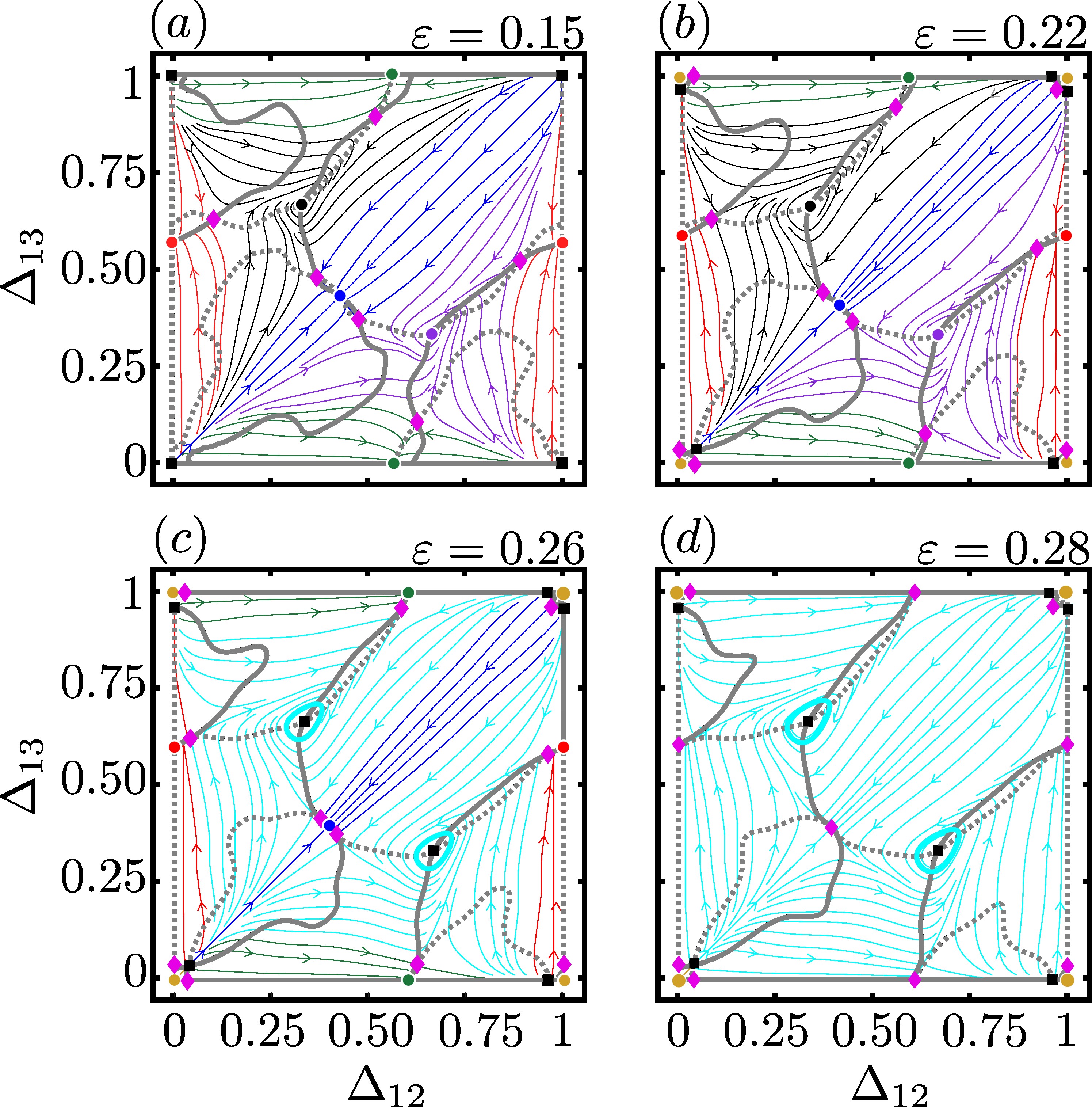}
	\caption{\label{var0.15}
		Polyrhythmic dynamics at \(\alpha=0.05\) for the parameter values \(A\)--\(D\) in Fig.~\ref{3d}(a). Three stable PMs coexist
		with two TWs in panel (a). SY become stable in panel~(b). Then TWs lose stability, and two invariant closed curves associated with LCs of
			the phase-difference system are present in panel~(c). The three PM states are no longer stable, while the two LCs persist in panel~(d).
		Arrows indicate the direction of the phase trajectories. The notation and color coding follow Figs.~\ref{compare} and~\ref{tau0.08}. Parameters: \(I=0.425\) and \(\varepsilon=0.15,0.22,0.26\) and \(0.28\) in panels~(a)--(d), resp.}
\end{figure}

To address this limitation, we employ standard phase-reduction theory. Under sufficiently weak coupling, the full network is approximated by a low-dimensional system of coupled phase oscillators. To first order, the reduced equations are determined by the infinitesimal phase-response curve (PRC) and the associated interaction function. The PRC quantifies the sensitivity of an oscillator's phase to a weak perturbation and can be determined numerically or experimentally. Once the limit cycle, PRC, and interaction function of the individual oscillator have been computed, the phase-reduction framework can be applied systematically to weakly coupled networks with different coupling architectures. The resulting phase-difference equations permit the direct calculation of phase-locked states, their stability, and the bifurcations that organize transitions among network rhythms.

\begin{figure}[b!]
	\centering
	\includegraphics[width=0.5\textwidth]{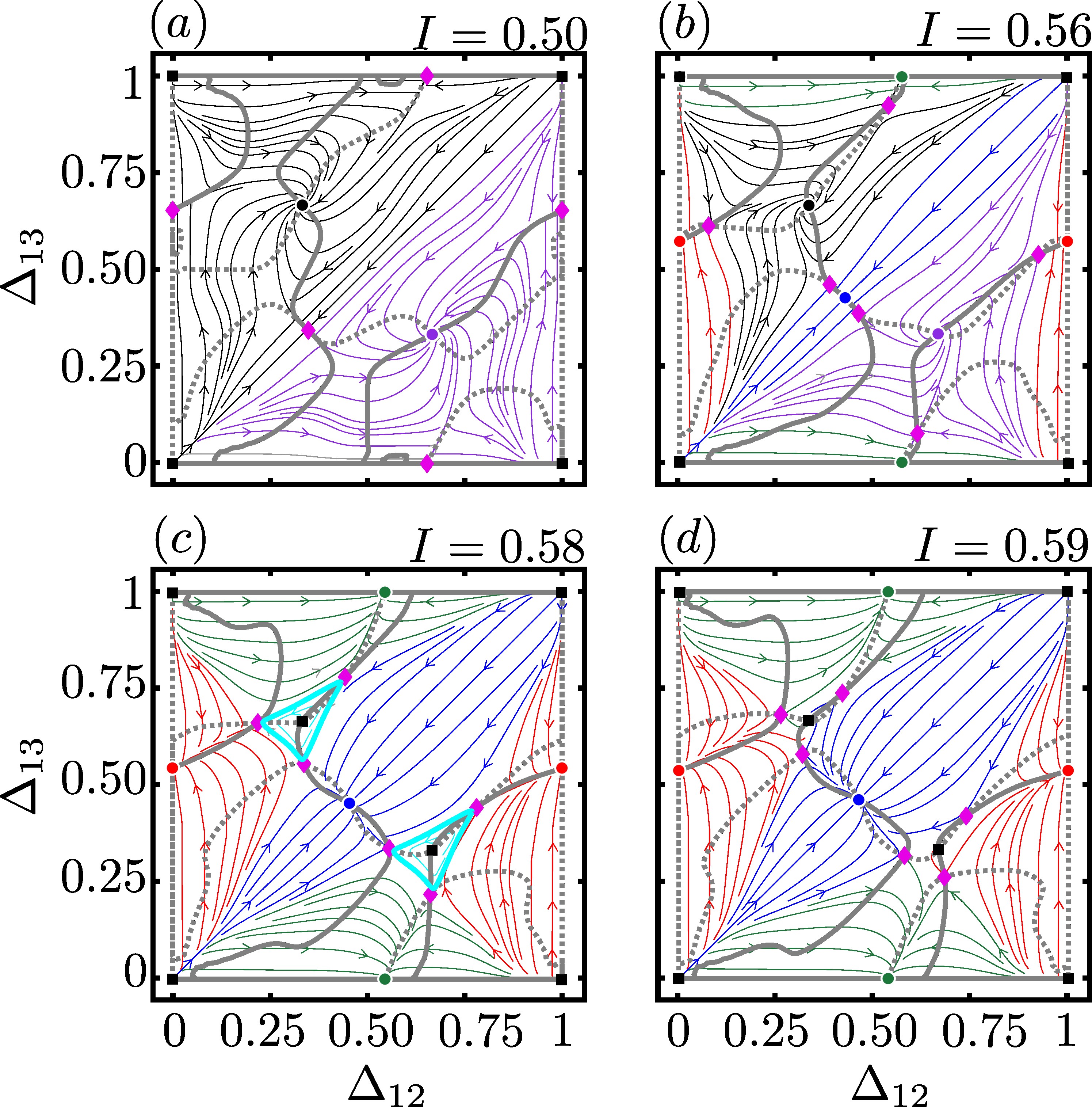}
	\caption{\label{var0.25}
		Polyrhythmic dynamics at \(\alpha=0.05\) for the parameter values $A'-D'$ in Fig.~\ref{3d}(a). The two TWs dominate the dynamics in panel~(a).
		As \(I\) increases, the three PMs appear in panel~(b).
		Then, the two TW FPs lose stability and two invariant closed curves appear in panel~(c). Finally, the PM retain the domination and these invariant closed curves disappear in panel~(d). The notation and color coding follow Figs.~\ref{compare} and~\ref{tau0.08}.  Parameters: \(\varepsilon=0.15\) and \(I=0.50,0.56,0.58\) and \(0.59\) in panels~(a)--(d), respectively.}
\end{figure}

We apply this framework in the weak-coupling regime \(g=0.0005\) to investigate how the intrinsic time-scale parameter \(\varepsilon\) and the applied current \(I\) jointly affect the rhythmic dynamics. These quantities are treated as bifurcation parameters in a two-parameter analysis of the phase model. 

Because changing either parameter modifies the intrinsic limit cycle, the corresponding PRC and interaction function are evaluated for each parameter pair, shown in Fig.~\ref{PRCs}. The lower dimensionality of the reduced system then permits a detailed classification of stable FPs and other phase-space structures across the \((\varepsilon,I)\) parameter plane.

The most pronounced changes in rhythmic organization occur for relatively short delays. We therefore focus on three representative cases occurring specifically at   \(\alpha=0.05\), \(0.10\), and \(0.20\), as shown in Fig.~\ref{3d}. The resulting diagrams reveal the global organization of the SY, PM, and TW states and their coexistence in the reduced phase space. Consistent with Ref.~\cite{zhu2026dynamical}, the rhythmic regimes exhibit an approximate reflection between the low- and high-current regions about an intermediate range of \(I\). This organization is not exactly symmetric, reflecting the dependence of the limit cycle and interaction function on the applied current.

\begin{figure}[t!]
	\centering
	\includegraphics[width=0.5\textwidth]{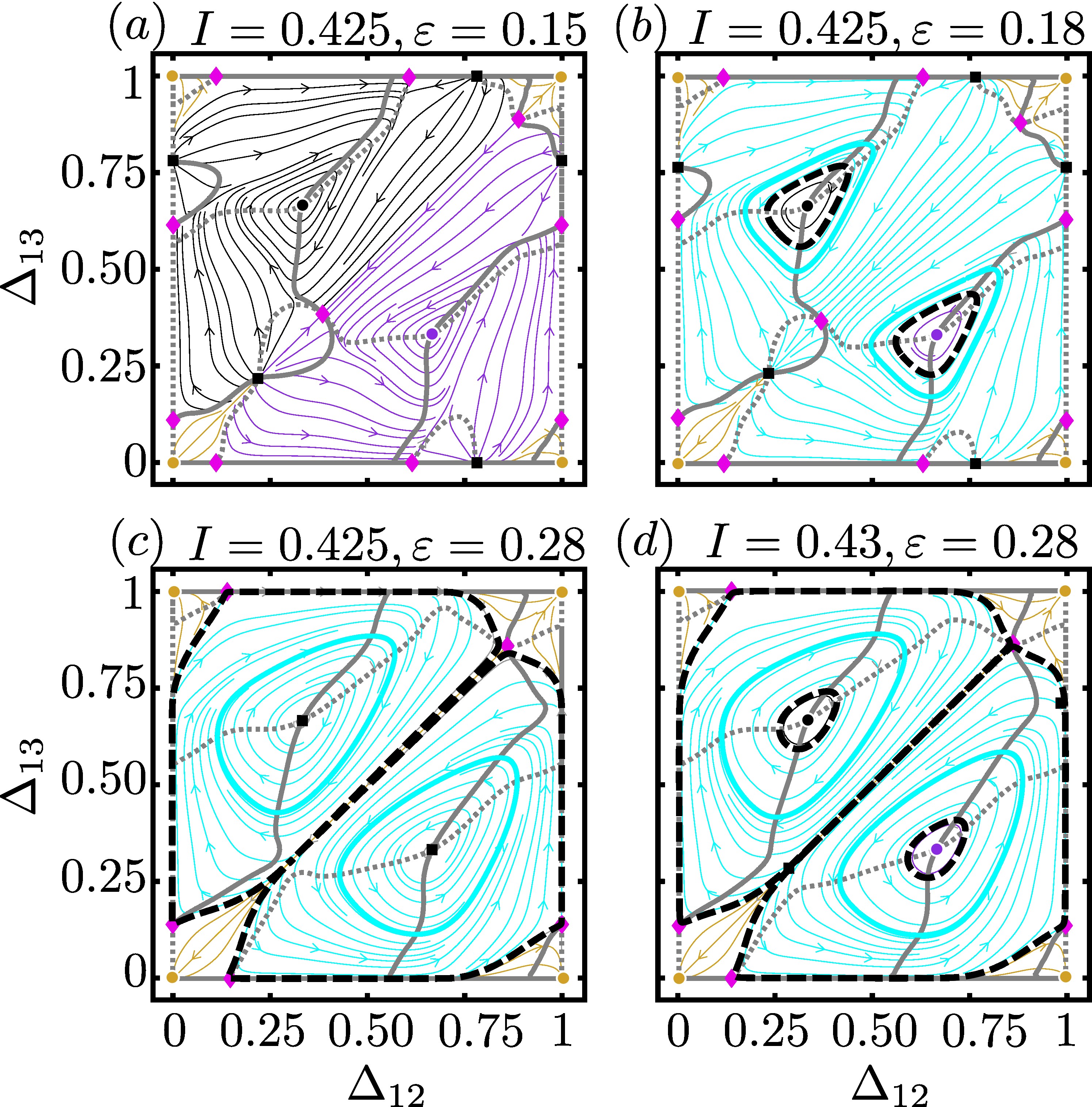}
	\caption{\label{tau0.1}
		Polyrhythmic dynamics at \(\alpha=0.10\) for the parameter values labelled \(A\)--\(D\) in Fig.~\ref{3d}(b). (a) Two stable TWs emerge, together with a SY with a small basin of attraction. (b) Two stable LCs (solid cyan curves) coexist with the SY attractor at the origin, and two unstable LCs (dashed black curves), nested around both TWs, separate attraction basins of the five coexisting attractors. (c) As \(\varepsilon\) increases, the TW FPs lose stability, the LCs dominate the dynamics, while the SY has a relatively small basin. (d) The two TWs regain stability, while the stable LCs persist and additional small unstable LCs emerge for larger \(I\)  values. The notation and color coding follow Figs.~\ref{compare} and~\ref{tau0.08}: the yellow dots stand for the stable SY at the origin, while the black square $\blacksquare$ and purple diamond $\blacklozenge$ symbols represent unstable FPs and saddles, and the cyan curves represent the stable invariant curves (ICs). Parameters: \(I=0.425\) and \(\varepsilon=0.15,0.18\) in panels~(a)--(b) and \(\varepsilon=0.28\) and \(I=0.425,0.43\) in panels~(c)--(d), resp.}
\end{figure}
\begin{figure}[t!]
	\centering
	\includegraphics[width=0.5\textwidth]{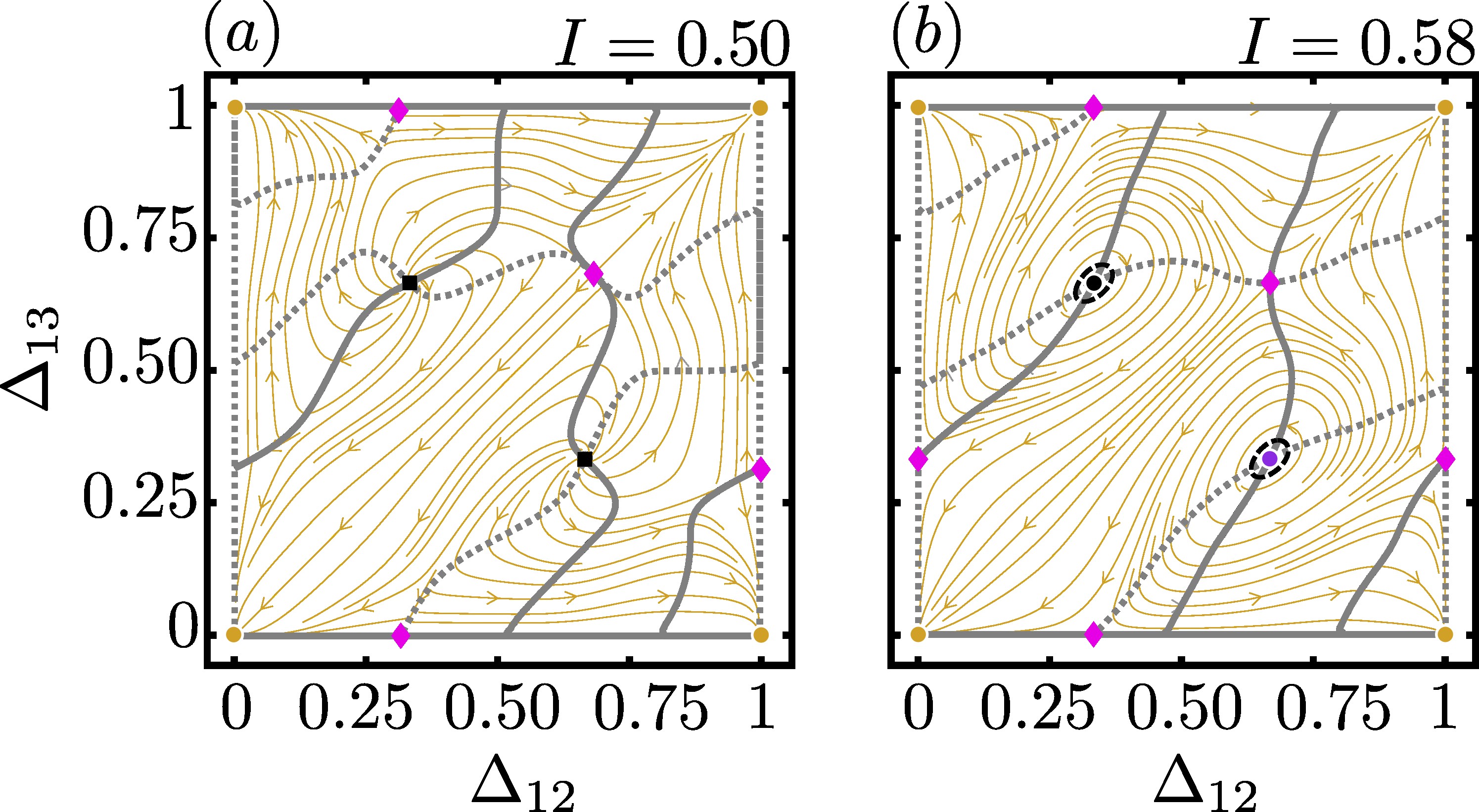}
	\caption{\label{tau0.2}	Polyrhythmic dynamics at \(\alpha=0.20\) for the parameter values \(A\) and \(B\) in the bifurcation diagram presented in Fig.~\ref{3d}(c). At \(I=0.50\), the SY state dominates the phase space in panel~(a). Increasing the current to \(I=0.58\) reveals two stable TWs coexisting with the stable SY state, while two unstable LCs emerge in panel~(b). The notation and color coding follow Figs.~\ref{compare} and~\ref{tau0.08}: the yellow dots stand for the stable SY at the origin, while the black square $\blacksquare$ and purple diamond $\blacklozenge$ symbols represent unstable FPs and saddles, and the cyan curves represent the stable invariant curves (ICs). Parameters: \(\varepsilon=0.20\) and \(I=0.50\) and \(0.58\) in panels~(a) and (b), resp.}
\end{figure}

At \(\alpha=0.05\), we examine the rhythmic dynamics along two parameter paths: varying \(\varepsilon\) at fixed \(I=0.425\) in Fig.~\ref{var0.15}, and varying \(I\) at fixed \(\varepsilon=0.15\) in Fig.~\ref{var0.25}.
In Fig.~\ref{var0.15}(a), three PMs and two TWs coexist with the colored trajectories showing their corresponding basins of
attraction.
As $\varepsilon$ increases, the two TWs change from stable in Fig.~\ref{var0.15}(b) to unstable foci and two invariant closed curves emerge in Fig.~\ref{var0.15}(c) through an Andronov--Hopf bifurcation. In the continuous phase-difference system, these invariant curves correspond to stable LCs.
Specifically, the corresponding eigenvalues change from $-0.0003\pm0.0094i$ to $0.0001\pm0.0119i$, indicating that the real part of a complex-conjugate pair crosses the imaginary axis at approximately $\varepsilon_H=0.255$.
The same Jacobian and eigenvalue analysis was also performed for the TWs in Figs.~\ref{var0.25}, ~\ref{tau0.1} and \ref{tau0.2}, confirming the occurrence of the corresponding Andronov--Hopf bifurcations. With a further increase in $\varepsilon$, the two saddles near each PM gradually approach and merge with it through a subcritical pitchfork bifurcation, causing the PM to lose stability. Consequently, the dynamics are dominated by stable LCs, while a stable SY state coexists with a small basin of attraction in Fig.~\ref{var0.15}(d).

For \(\varepsilon=0.15\), two stable TWs dominate the dynamics at \(I=0.50\) [Fig.~\ref{var0.25}(a)]. Then three PMs appear and coexist with two TWs in Fig.~\ref{var0.25}(b). When \(I\) is increased to \(0.58\), the TW FPs lose stability while two LCs emerge through Andronov--Hopf bifurcations around their locations in the return map in Fig.~\ref{var0.25}(c). Finally, the LCs disappear and PM states become the dominate rhythm in Fig.~\ref{var0.25}(d).

Figure~\ref{tau0.1} illustrates transitions among multiple coexisting attractors at the fixed delay \(\alpha=0.10\). At \(\varepsilon=0.15\), two stable TWs dominate most of the phase space, while a smaller set of initial conditions converges to the stable SY state [Fig.~\ref{tau0.1}(a)]. As \(\varepsilon\) increases to \(0.18\), both two stable LCs emerge, thus coexisting with the three persisting attractors: two TWs and the SY state, as documented in Fig.~\ref{tau0.1}(b). 

With a further increase in $\varepsilon$, the small unstable LCs gradually shrink toward the stable TW equilibria and eventually disappear as the TWs lose stability through subcritical Andronov--Hopf bifurcations, while the two stable LCs persist. Meanwhile, the unstable fixed points approach and eventually collide with nearby saddles, leading to their annihilation, while larger unstable LCs emerge via homoclinic saddle-node bifurcation in Fig.~\ref{tau0.1}(c).
For fixed $\varepsilon=0.28$, as $I$ increases to 0.43, two TWs undergo Andronov–Hopf bifurcations and regain stability, accompanied by the emergence of two small unstable LCs surrounding them. The larger unstable LCs and the two stable LCs persist, as shown in panel (d). Now, the phase-difference system exhibits five coexisting attractors: two stable TWs, two stable LCs, and one stable SY state. The small unstable LCs separate the basins of attraction of the stable TWs and the surrounding stable LCs, while the larger unstable LCs delimit the basin of attraction of the SY state.

Finally, we examine the dynamics at \(\alpha=0.20\). For \(I=0.50\), the SY state dominates the phase space, see Fig.~\ref{tau0.2}(a). When the current is increased to \(I=0.58\), both repelling TWs regain stability through Andronov--Hopf bifurcations, accompanied by the emergence of two small unstable LCs surrounding them, and now coexist with the synchronous rhythm [Fig.~\ref{tau0.2}(b)].

\section{Conclusion}
We investigated the dynamics and transformations of delayed three-cell inhibitory CPG networks by combining return-map analysis of the full generalized Fitzhugh--Nagumo system with a corresponding phase-difference model. The two approaches provide complementary descriptions of the network dynamics. Return maps retain a direct connection to the full system and reveal the global organization of rhythmic trajectories and their basins of attraction. Within the weak-coupling approximation, the phase-difference model provides a more computationally efficient and quantitatively precise description of the relative-phase dynamics. It permits direct localization of phase-locked states through nullcline intersections, determination of their stability from the linearized system, and identification of stable limit cycles corresponding to invariant closed curves observed in the return maps.

The phase model is particularly useful when an attractor has a very small basin of attraction. Such states can be difficult to detect through direct forward simulations because only a small set of initial conditions converges to them. In contrast, the corresponding equilibria of the phase-difference system can be located and classified independently of their basin size. This distinction allowed us to identify the synchronous FP and determine its stability even in parameter regimes where its basin is extremely small.

Using the return-map framework, we first examined the joint effects of the normalized delay \(\alpha\) and the reciprocal coupling strength \(g_{13/31}\), varied relative to the remaining inhibitory synapses. The resulting parameter scans revealed a rich multistable landscape containing synchronous, pacemaker, traveling-wave, phase-slipping, and limit-cycle dynamics, together with numerous coexistence regimes. Phase-slipping dynamics persist across a broad region of parameter space, while variations in coupling asymmetry and delay organize transitions among synchronous, traveling-wave, and pacemaker rhythms.

We then used the phase-reduced model to perform two-parameter scans in the intrinsic time-scale parameter \(\varepsilon\) and the applied current \(I\) for relatively short delays. The reduced system reproduced the global organization of the principal rhythms and revealed an approximate reflection between the low- and high-current regimes. This reflection is accompanied by an approximate reversal of the PRC profiles between low and high $I$-value. Because the reduced equations evolve only the phase differences, they permit denser parameter sweeps at substantially lower computational cost than simulations of the full delayed network. This efficiency makes it possible to resolve narrow stability regions, small-basin attractors, and bifurcation structures that can be difficult to identify using return maps alone.

Together, return-map analysis and phase reduction provide a unified framework for characterizing multistability, bifurcations, and transitions among rhythms in delayed inhibitory CPG networks. Their close agreement validates the phase-reduced description in the weak-coupling regime. The phase model should not, however, be viewed as a replacement for the full-system analysis: return maps remain essential when stronger coupling or delay-induced waveform changes produce amplitude effects that lie outside the assumptions of first-order phase reduction.

Future work should extend this framework to larger and less symmetric network architectures to determine how the number, stability, and organization of coexisting rhythms scale with network size and connectivity. A more detailed analysis of the local and global bifurcations responsible for the creation, stability changes, and disappearance of invariant curves would also clarify the mechanisms underlying transitions between phase locking, modulated rhythms, and persistent phase slipping.

\section*{ACKNOWLEDGMENT}
X.Q.'s work partly funded by the National Natural Science Foundation of China (Grant Nos. 11872183) and China Scholarship Council (Grant Nos. 202506150077).
A.L.S's research was partially funded by a National Science Foundation (NSF) (Award No.
DMS-2407999).

\section*{DATA AVAILABILITY}
Some of the code used for the numerical simulations is publicly available at [\url{https://github.com/xinxin77-cloud/FHN3_phase_difference}]. 

\bibliographystyle{apsrev4-2}
\bibliography{our_refs} 

@article{marder2001central,
  author = {Marder, Eve and Bucher, Dirk},
  journal = {Current biology},
  number = {23},
  pages = {R986--R996},
  publisher = {Elsevier},
  title = {Central pattern generators and the control of rhythmic movements},
  volume = {11},
  year = {2001}
}

@article{yuste2005cortex,
  author = {Yuste, Rafael and MacLean, Jason N and Smith, Jeffrey and Lansner, Anders},
  journal = {Nature Reviews Neuroscience},
  number = {6},
  pages = {477--483},
  publisher = {Nature Publishing Group UK London},
  title = {The cortex as a central pattern generator},
  volume = {6},
  year = {2005}
}

@article{kiehn2003physiological,
  author = {Kiehn, Ole and Butt, Simon JB},
  journal = {Progress in neurobiology},
  number = {4},
  pages = {347--361},
  publisher = {Elsevier},
  title = {Physiological, anatomical and genetic identification of CPG neurons in the developing mammalian spinal cord},
  volume = {70},
  year = {2003}
}

@article{qie2025synaptic,
  author = {Qie, Xinxin and Zang, Jie and Liu, Shenquan and Shilnikov, Andrey L},
  journal = {Chaos: An Interdisciplinary Journal of Nonlinear Science},
  number = {4},
  publisher = {AIP Publishing},
  title = {Synaptic delays shape dynamics and function in multimodal neural motifs},
  volume = {35},
  year = {2025}
}

@article{shafiei2020time,
  author = {Shafiei, Mohadeseh and Jafari, Sajad and Parastesh, Fatemeh and Ozer, Mahmut and Kapitaniak, Tomasz and Perc, Matja{\v{z}}},
  journal = {Communications in Nonlinear Science and Numerical Simulation},
  pages = {105175},
  publisher = {Elsevier},
  title = {Time delayed chemical synapses and synchronization in multilayer neuronal networks with ephaptic inter-layer coupling},
  volume = {84},
  year = {2020}
}

@inproceedings{yu2018effect,
  author = {Yu, Minxing and Wang, Jiang and Yu, Haitao and Qin, Yingmei},
  booktitle = {2018 37th Chinese Control Conference (CCC)},
  organization = {IEEE},
  pages = {8183--8187},
  title = {Effect of Time Delay on the variability of output pattern of CPG network},
  year = {2018}
}

@incollection{venugopal2022connectionist,
  author = {Venugopal, Sharmila},
  booktitle = {Encyclopedia of Computational Neuroscience},
  pages = {957--959},
  publisher = {Springer},
  title = {Connectionist Models of CPG Networks},
  year = {2022}
}

@article{song2026coupling,
  author = {Song, Caihong and Song, Zigen and Kong, Qingzhao and Sun, Xiuting and Xu, Jian},
  journal = {Nonlinear Dynamics},
  number = {4},
  pages = {268},
  publisher = {Springer},
  title = {Coupling delays regulate synchrony patterns to replicate the function of CPG networks: C. Song et al.},
  volume = {114},
  year = {2026}
}

@article{campbell2012phase,
  author = {Campbell, Sue Ann and Kobelevskiy, Ilya},
  journal = {Discrete Contin. Dyn. Syst},
  pages = {2653--2673},
  title = {Phase models and oscillators with time delayed coupling},
  volume = {8},
  year = {2012}
}

@article{wojcik2014key,
  author = {Wojcik, Jeremy and Schwabedal, Justus and Clewley, Robert and Shilnikov, Andrey L},
  journal = {PloS one},
  number = {4},
  pages = {e92918},
  publisher = {Public Library of Science San Francisco, USA},
  title = {Key bifurcations of bursting polyrhythms in 3-cell central pattern generators},
  volume = {9},
  year = {2014}
}

@article{kelley20202,
  author = {Kelley, Aaron and Shilnikov, Andrey L},
  journal = {Frontiers in Applied Mathematics and Statistics},
  pages = {588904},
  publisher = {Frontiers Media SA},
  title = {2 $\theta$-burster for rhythm-generating circuits},
  volume = {6},
  year = {2020}
}

@article{collens2020dynamics,
  author = {Collens, J and Pusuluri, K and Kelley, A and Knapper, D and Xing, T and Basodi, S and Alacam, Deniz and Shilnikov, Andrey L},
  journal = {Chaos: An Interdisciplinary Journal of Nonlinear Science},
  number = {7},
  publisher = {AIP Publishing},
  title = {Dynamics and bifurcations in multistable 3-cell neural networks},
  volume = {30},
  year = {2020}
}

@article{wojcik2011order,
  author = {Wojcik, Jeremy and Clewley, Robert and Shilnikov, Andrey L},
  journal = {Physical Review E—Statistical, Nonlinear, and Soft Matter Physics},
  number = {5},
  pages = {056209},
  publisher = {APS},
  title = {Order parameter for bursting polyrhythms in multifunctional central pattern generators},
  volume = {83},
  year = {2011}
}

@book{ermentrout2010mathematical,
  author = {Ermentrout, Bard and Terman, David M},
  publisher = {Springer},
  title = {Mathematical foundations of neuroscience},
  volume = {35},
  year = {2010}
}

@inproceedings{kuramoto2005self,
  author = {Kuramoto, Yoshiki},
  booktitle = {International symposium on mathematical problems in theoretical physics: January 23--29, 1975, Kyoto university, Kyoto/Japan},
  organization = {Springer},
  pages = {420--422},
  title = {Self-entrainment of a population of coupled non-linear oscillators},
  year = {2005}
}

@article{scully2024pairing,
  author = {Scully, James and Bourahmah, Jassem and Bloom, David and Shilnikov, Andrey L},
  journal = {Frontiers in Network Physiology},
  pages = {1397151},
  publisher = {Frontiers Media SA},
  title = {Pairing cellular and synaptic dynamics into building blocks of rhythmic neural circuits. A tutorial},
  volume = {4},
  year = {2024}
}

@article{zhu2026dynamical,
  author = {Zhu, Qinghua and Zheng, Song and Han, Fang},
  journal = {Nonlinear Dynamics},
  number = {8},
  pages = {593},
  publisher = {Springer},
  title = {Dynamical analysis of multistable 3-cell CPG based on phase response curve and interaction function: Q. Zhu et al.},
  volume = {114},
  year = {2026}
}

@article{schuster1989mutual,
  author = {Schuster, Heinz Georg and Wagner, Peter},
  journal = {Progress of Theoretical Physics},
  number = {5},
  pages = {939--945},
  publisher = {Oxford University Press},
  title = {Mutual entrainment of two limit cycle oscillators with time delayed coupling},
  volume = {81},
  year = {1989}
}

@article{yeung1999time,
  author = {Yeung, MK Stephen and Strogatz, Steven H},
  journal = {Physical review letters},
  number = {3},
  pages = {648},
  publisher = {APS},
  title = {Time delay in the Kuramoto model of coupled oscillators},
  volume = {82},
  year = {1999}
}

@article{bressloff1997phase,
  author = {Bressloff, Paul C and Coombes, S},
  journal = {Physica D: Nonlinear Phenomena},
  number = {3-4},
  pages = {232--254},
  publisher = {Elsevier},
  title = {Travelling waves in chains of pulse-coupled integrate-and-fire oscillators with distributed delays},
  volume = {130},
  year = {1999}
}

@article{earl2003synchronization,
  author = {Earl, Matthew G and Strogatz, Steven H},
  journal = {Physical Review E},
  number = {3},
  pages = {036204},
  publisher = {APS},
  title = {Synchronization in oscillator networks with delayed coupling: A stability criterion},
  volume = {67},
  year = {2003}
}

@article{crook1997role,
  author = {Crook, Sharon M. and Ermentrout, G Bard and Vanier, Michael C. and Bower, James M.},
  journal = {Journal of computational neuroscience},
  number = {2},
  pages = {161--172},
  publisher = {Springer},
  title = {The role of axonal delay in the synchronization of networks of coupled cortical oscillators},
  volume = {4},
  year = {1997}
}

@article{novivcenko2012phase,
  author = {Novi{\v{c}}enko, V and Pyragas, K},
  journal = {Physica D: Nonlinear Phenomena},
  number = {12},
  pages = {1090--1098},
  publisher = {Elsevier},
  title = {Phase reduction of weakly perturbed limit cycle oscillations in time-delay systems},
  volume = {241},
  year = {2012}
}

@article{kotani2012adjoint,
  author = {Kotani, Kiyoshi and Yamaguchi, Ikuhiro and Ogawa, Yutaro and Jimbo, Yasuhiko and Nakao, Hiroya and Ermentrout, G Bard},
  journal = {Physical Review Letters},
  number = {4},
  pages = {044101},
  publisher = {APS},
  title = {Adjoint method provides phase response functions for delay-induced oscillations},
  volume = {109},
  year = {2012}
}

@article{kurebayashi2022phase,
  author = {Kurebayashi, Wataru and Yamamoto, Takuto and Shirasaka, Sho and Nakao, Hiroya},
  journal = {Physical Review Research},
  number = {4},
  pages = {043176},
  publisher = {APS},
  title = {Phase reduction of strongly coupled limit-cycle oscillators},
  volume = {4},
  year = {2022}
}

@article{pusuluri2020computational,
  author = {Pusuluri, Krishna and Basodi, Sunitha and Shilnikov, Andrey L},
  journal = {Communications in Nonlinear Science and Numerical Simulation},
  pages = {105139},
  publisher = {Elsevier},
  title = {Computational exposition of multistable rhythms in 4-cell neural circuits},
  volume = {83},
  year = {2020}
}

@article{schwabedal2016qualitative,
  author = {Schwabedal, Justus TC and Knapper, Drake E and Shilnikov, Andrey L},
  journal = {Nonlinearity},
  number = {12},
  pages = {3647--3676},
  publisher = {IOP Publishing},
  title = {Qualitative and quantitative stability analysis of penta-rhythmic circuits},
  volume = {29},
  year = {2016}
}

@book{ermentrout02xppaut,
	address = {Philadelphia},
	author = {Ermentrout, G.B.},
	publisher = {SIAM Books},
	title = {Simulating, {A}nalyzing, and {A}nimating {D}ynamical {S}ystems: {A} {G}uide to {XPPAUT} for {R}esearchers and {S}tudents},
    doi={10.1137/1.9780898718195.fm},
	year = {2002}
}

@article{Izhikevich1998,
   author = {Eugene M. Izhikevich},
   doi = {10.1103/PhysRevE.58.905},
   issn = {1063-651X},
   issue = {1},
   journal = {Physical Review E},
   month = {7},
   pages = {905-908},
   title = {Phase models with explicit time delays},
   volume = {58},
   year = {1998}
}

\end{document}